\documentclass[a4paper]{article}
\usepackage{amscd,amsmath,amsthm,amssymb}
\begin{document}

\newtheorem{thm}{Theorem}[section]
\newtheorem{lem}[thm]{Lemma}
\newtheorem{prop}[thm]{Proposition}
\newtheorem{cor}[thm]{Corollary}
\newtheorem{bigthm}{Theorem}
\renewcommand{\thebigthm}{\Alph{bigthm}}
\newtheorem{bigcor}[bigthm]{Corollary}
\newtheorem{exam}[thm]{Example}

\newenvironment{pf}{\paragraph{Proof}}{\par\smallskip}
\newenvironment{pfof}[1]{\paragraph{Proof of #1}}{\par\smallskip}
\newenvironment{example}{\begin{exam}\em}{\end{exam}}

\newcommand{\dual}{^{\vee}}
\newcommand{\contr}{\vdash}
\newcommand{\de}{\partial}
\newcommand{\debar}{{\overline{\partial}}}
\newcommand{\desude}[2]{{\ds\frac{\de #1}{\de #2}}}
\newcommand{\cat}[1]{\begin{bf}#1\end{bf}}
\newcommand{\mapor}[1]{\smash{\mathop{\longrightarrow}\limits^{#1}}}
\newcommand{\Art}{\cat{Art}}
\newcommand{\Set}{\cat{Set}}
\renewcommand{\bar}{\overline}
\newcommand{\h}[1]{\widehat{#1}}
\renewcommand{\t}[1]{\widetilde{#1}}
\newcommand{\ds}{\displaystyle}
\newcommand{\oo}{\infty}
\newcommand{\sA}{{\mathcal A}}
\newcommand{\Oh}{{\mathcal O}}
\newcommand{\sH}{{\mathcal H}}
\newcommand{\sF}{{\mathcal F}}
\newcommand{\sM}{{\mathcal M}}
\newcommand{\sG}{{\mathcal G}}
\newcommand{\sX}{{\mathcal X}}
\newcommand{\sB}{{\mathcal B}}
\newcommand{\sY}{{\mathcal Y}}

\newcommand{\C}{\mathbb C}
\newcommand{\Z}{\mathbb Z}
\newcommand{\MC}{\operatorname{MC}}
\newcommand{\Def}{\operatorname{Def}}
\newcommand{\Hom}{\operatorname{Hom}}
\newcommand{\Der}{\operatorname{Der}}
\newcommand{\somdir}[2]{\hbox{$\mathrel
{\smash{\mathop{\mathop \bigoplus\limits _{#1}}
\limits^{#2}}}$}}
\newcommand{\tensor}[2]{\hbox{$\mathrel
{\smash{\mathop{\mathop \bigotimes\limits _{#1}}
^{#2}}}$}}
\newcommand{\symm}[2]{\hbox{$\mathrel
{\smash{\mathop{\mathop \bigodot\limits _{#1}}
^{#2}}}$}}

\title{Cohomological constraint to deformations of compact 
K\"ahler manifolds}
\author{Marco Manetti\thanks{Partially 
supported by Italian MURST-PRIN 
'Spazi di moduli e teoria delle rappresentazioni'. Member of GNSAGA 
of CNR.}\\
Universit\`a di Roma ``La Sapienza'', Italy}
\date{}

\maketitle

\begin{center}
{\it Dedicated to the memory of Fabio Bardelli}
\end{center}
\bigskip

\begin{abstract}
We prove that for every compact K\"ahler manifold $X$ the
cup product 
\[H^*(X,T_X)\otimes H^*(X,\Omega_X^*)\to 
H^*(X,\Omega_X^{*-1})\]
can be lifted to   
an $L_{\oo}$-morphism from the Kodaira-Spencer differential 
graded Lie algebra to the suspension of the 
space of linear endomorphisms of the singular cohomology of $X$.
As a consequence we get an algebraic proof of the principle 
``obstructions to deformations of compact K\"{a}hler manifolds 
annihilate ambient cohomology''.
\end{abstract}

\section*{Introduction}

In this paper we give an algebraic proof 
of the principle ``obstructions to deformations of compact K\"{a}hler manifolds 
annihilate ambient cohomology'' recently proved, in a different way, 
by Herb Clemens~\cite{clemens} and Ziv Ran~\cite{ran}.\\ 
Let $X$ be a fixed compact K\"ahler manifold of dimension $n$ and 
consider the graded vector space 
$M_{X}=\Hom_{\C}^{*}(H^{*}(X,\C),H^{*}(X,\C))$  
of linear endomorphisms of the singular cohomology of $X$. The 
Hodge decomposition gives natural isomorphisms 
\[ M_{X}=\somdir{i}{}M^{i}_{X},\qquad 
M^{i}_{X}=\somdir{r+s=p+q+i}{}\Hom_{\C}
(H^{p}(\Omega^{q}_{X}),H^{r}(\Omega^{s}_{X}))\]
and the composition of the cup product and the contraction 
operator $T_{X}\otimes\Omega^{p}_{X}\mapor{\contr}\Omega^{p-1}_{X}$ 
gives natural linear maps 
\[ \theta_{p}\colon H^{p}(X,T_{X})\to 
\somdir{r,s}{}\Hom_{\C}^{*}
(H^{r}(\Omega^{s}_{X}),H^{r+p}(\Omega^{s-1}_{X}))
\subset M[-1]_{X}^{p}=M^{p-1}_{X}.\]
The Dolbeaut's complex of the holomorphic tangent bundle $T_{X}$ 
\[ KS_{X}=\somdir{p}{}KS^{p}_{X},\qquad 
KS^{p}_{X}=\Gamma(X,\sA^{0,p}(T_{X}))\]
has a natural structure of differential graded Lie algebra (DGLA), 
\cite{Catacime}, \cite{GoldMil2}, \cite[3.4.1]{K},
called the Kodaira-Spencer algebra of $X$. By Dolbeaut's theorem 
$H^{*}(KS_{X})=H^{*}(X,T_{X})$ and then the maps $\theta_{i}$ give
a morphism of graded 
vector spaces $\theta\colon H^{*}(KS_{X})\to M[-1]_{X}$.
This morphism is  generally nontrivial: consider for instance a 
Calabi-Yau manifold where the map $\theta_{p}$ induces an
isomorphism
$H^{p}(X,T_{X})=\Hom_{\C}(H^{0}(\Omega^{n}_{X}),H^{p}(\Omega^{n-1}_{X}))$.

\begin{bigthm}\label{thmA} In the above notation, 
consider $M[-1]_{X}$ as a differential graded Lie algebra with trivial 
differential and trivial bracket.\\ 
Every choice of a K\"ahler metric 
on $X$ induces a canonical lifting of $\theta$ to an $L_{\oo}$-morphism 
from $KS_{X}$ to $M[-1]_{X}$.				
\end{bigthm}

The above theorem, together some standard and purely formal results 
in Schlessinger's theory, gives  immediate applications to the study of 
deformations of $X$.  
In fact the deformations of $X$ are governed by the Kodaira-Spencer 
differential 
graded Lie algebra $KS_{X}$ and every $L_{\oo}$-morphism between DGLAs
induces a natural transformation between the associated deformation
functors.  The triviality of the DGLA structure on $M[-1]_{X}$ allows to 
prove easily the following:

\begin{bigcor}\label{CorB} Let $f\colon\sY\to \sB$ be the semiuniversal 
deformation  of  a compact K\"ahler manifold $Y$ and let 
$X\mapor{\pi}Y$ be a finite unramified covering. 
For every $p\ge 0$ denote by $\alpha_p$ the composite linear map 
\[\alpha_p\colon H^{p}(Y,T_{Y})\mapor{\pi^{*}}H^{p}(X,T_{X})
	\mapor{\theta_{p}}\somdir{r,s}{}
\Hom_{\C}(H^{r}(\Omega^{s}_{X}),H^{r+p}(\Omega^{s-1}_{X})).\]
Then: 
\begin{enumerate}
	\item If $\alpha_1$
is injective then $f\colon\sY\to \sB$ is universal.
	\item There exists a morphisms of complex analytic singularities 
	$q\colon (H^{1}(Y,T_{Y}),0)\to (\ker\alpha_{2},0)$
	such that $\sB$ is isomorphic to $q^{-1}(0)$.
In particular if  $\alpha_{2}$ is injective  then $\sB$ is smooth. 
\end{enumerate}
\end{bigcor}

As an example, if $Y$ is a projective manifold with torsion canonical 
bundle and $\pi\colon X\to Y$ is the canonical covering, then all the 
maps $\alpha_{p}$ are injective.\\
 
Probably the main interesting aspect of  Theorem~\ref{thmA} is that 
it gives a concrete construction of a 
morphism whose existence is predicted by the general philosophy of 
extended deformation theory.\\
Roughly speaking, to every deformation problem over a field of 
characteristic 0, it is associated a differential  graded Lie algebra 
$L$, unique up to quasiisomorphism,
and a formal pointed quasismooth dg-manifold $\sM$ quasiisomorphic to 
$L$ as $L_{\oo}$-algebra. The differential graded Lie algebra $L$ 
governs the deformation problem via the solutions Maurer-Cartan modulo 
gauge action and the truncation in degree 0 of $\sM$ is the classical    
moduli space (cf. \cite{EDF}, 
Section 2 of \cite{BreSo} and references therein).\\
Moreover, according to this general philosophy, every natural 
morphism between moduli spaces (e.g. the period map from deformations 
of a compact K\"ahler manifold to deformations of its Hodge decomposition) 
should extend to a morphism of their  extended moduli spaces and 
therefore induces an $L_{\oo}$-morphism between the associated 
differential graded Lie algebras.\\

The author thanks A. Canonaco for his useful 
help in the preparation of the paper.

\subsection*{Notation}
For every holomorphic vector bundle $E$ on a complex manifold 
we denote by $\sA^{p,q}(E)$ the sheaf 
of differential $(p,q)$-forms with coefficients in $E$.\\
For every vector space $V$ and every linear functional $\alpha\colon 
V\to \C$ we denote by $\alpha\contr\colon\bigwedge^{k}V\to 
\bigwedge^{k-1}V$ the contraction operator
\[\alpha\contr(v_{1}\wedge\ldots\wedge v_{k})=
\sum_{i=1}^{k}(-1)^{i-1}\alpha(v_{i})v_{1}\wedge\ldots
\wedge\h{v_{i}}\wedge\ldots\wedge v_{k}.\]
We point out for later use that $\alpha\contr$ is a derivation of 
degree $-1$ of the graded algebra 
$(\bigwedge^{*}V,\wedge)$.\\ 
We denote by $\Sigma_{m}$ the symmetric group of permutations of the 
set $\{1,2,\ldots,m\}$ and, for every $0\le p\le m$ by 
$S(p,m-p)\subset\Sigma_{m}$ the set of unshuffles of type $(p,m-p)$. 
By definition $\sigma\in S(p,m-p)$ if and only if 
$\sigma_{1}<\sigma_{2}<\ldots<\sigma_{p}$ and 
$\sigma_{p+1}<\sigma_{p+2}<\ldots<\sigma_{m}$.\\

\bigskip
\section{$L_{\oo}$-morphisms}

Let $V=\oplus V^{i}$ be a $\Z$-graded vector space, for every 
integer $n$ we denote by 
$V[n]=\oplus V[n]^{i}$ the graded vector space where $V[n]^{i}=V^{n+i}$. 
The space $V[-1]$ is also called the suspension of $V$ and $V[1]$ the 
unsuspension.\\ 
The graded $m$-th symmetric power of $V$ is denoted by $\symm{}{m}V$. 
If $\sigma\in \Sigma_{m}$ and $a_{1},\ldots,a_{m}\in V$ are 
homogeneous elements, the Koszul sign 
$\epsilon(V,\sigma;a_{1},\ldots,a_{m})=\pm 1$ is defined by the rule
\[ a_{\sigma_1}\odot\ldots\odot a_{\sigma_m}=
\epsilon(V,\sigma;a_{1},\ldots,a_{m})a_{1}\odot\ldots\odot a_{m}\in 
\symm{}{m}V.\]
For simplicity of notation we write $\epsilon(V,\sigma)$ when the 
elements $a_{1},\ldots,a_{m}$ are clear from the context. If $a\in V$ 
is homogeneous we denote by $\deg(a,V)$ its degree; we also write 
$\deg(a,V)=\bar{a}$ when there is no ambiguity about $V$. Note that 
$\deg(a,V[n])=\deg(a,V)-n$.\\
We denote by $C(V)$ the reduced graded symmetric coalgebra generated by 
$V[1]$; more precisely it is the graded vector space 
\[ C(V)=\bar{S}(V[1])=\somdir{m=1}{\oo}\symm{}{m}(V[1])\]
endowed with the coproduct $\Delta\colon C(V)\to C(V)\otimes C(V)$, 
$\Delta(a)=0$ for every $a\in V[1]$ and 
\[ \Delta(a_{1}\odot\ldots\odot a_{m})=\sum_{r=1}^{m-1}
\sum_{\sigma\in S(r,m-r)}\!\!\!\!\!\!\epsilon(V[1],\sigma)
(a_{\sigma_1}\odot\ldots\odot a_{\sigma_r})\otimes
(a_{\sigma_{r+1}}\odot\ldots\odot a_{\sigma_m})\]
for every $a_{1},\ldots,a_{m}\in V[1]$, $m\ge 2$.\\
Assume now that $V$ has a structure of differential graded Lie 
algebra with differential $d$ and bracket $[\,\,  ,\,]$, then the linear map 
\[ Q\colon \symm{}{2}(V[1])\to V[1],\qquad 
Q(a\odot b)=(-1)^{\deg(a,V[1])}[a,b]\]
has degree 1 and the map $\delta\colon C(V)\to C(V)$ defined by 
\begin{equation}\label{codiff}
\begin{split}
\delta(a_{1}\odot\ldots\odot a_{m})=&\sum_{\sigma\in S(1,m-1)}
\epsilon(V[1],\sigma; a_{1},\ldots,a_{m})da_{\sigma_1}\odot 
a_{\sigma_2}\odot\ldots\odot
a_{\sigma_m}+\\
&+\!\!\!\!\!\!\sum_{\sigma\in S(2,m-2)}
\!\!\!\!\!\!\!\!\epsilon(V[1],\sigma; a_{1},\ldots,a_{m})Q(a_{\sigma_1}
\odot a_{\sigma_2})\odot a_{\sigma_3}\odot\ldots\odot
a_{\sigma_m}
\end{split}
\end{equation}
is a codifferential of degree 1 on the 
coalgebra $C(V)$. The differential graded coalgebra $(C(V),\delta)$ is 
called the $L_{\oo}$-algebra associated to the DGLA $(V,d, [~,~])$.\\
By definition, an $L_{\oo}$-morphism between two DGLA $V,V'$ is a 
morphism of differential graded coalgebras $\Theta\colon 
(C(V),\delta)\to (C(V'),\delta')$.\\
It is easy to check that if $f\colon V\to V'$ is a morphism of 
differential graded Lie algebras then the linear map 
\[ (C(V),\delta)\to (C(V'),\delta'),\qquad 
a_{1}\odot\ldots\odot a_{m}\to f(a_{1})\odot\ldots\odot f(a_{m})\]
is an $L_{\oo}$-morphism.
We refer to 
\cite{K}, \cite{LadaMarkl}, \cite{LadaStas}, \cite{SchSta} 
for the general theory of 
$L_{\oo}$-morphisms.\\ 
In this paper we are interested only in the 
particular and simple 
case when $V'$ has trivial differential and trivial bracket: 
under these assumption $\delta'=0$ and there exists a bijection 
between the set of $L_{\oo}$-morphism $\Theta\colon (C(V),\delta)\to
(C(V'),0)$ and the set of morphisms of graded vector spaces 
$F\colon C(V)\to V'[1]$ such that $F\circ \delta=0$. The bijection is 
given by the formulas
\[ F=p_{1}\circ\Theta,\qquad p_{1}\colon C(V')\to \symm{}{1}V'[1]=V'[1]
\hbox{~~ the projection}\]
\[\Theta=\sum_{m=1}^{\oo}\frac{1}{m!}F^{\odot m}\circ\Delta^{m-1}_{C(V)}
\colon C(V)\to C(V')\] 
where $F^{\odot m}$ is the composition of $F^{\otimes m}\colon 
\tensor{}{m}C(V)\to\tensor{}{m}(V'[1])$ with the projection onto 
the symmetric product  
$\tensor{}{m}(V'[1])\to\symm{}{m}(V'[1])$.\\
Let $F_{1}\colon V[1]\to V'[1]$ the composition 
of $F$ with the inclusion $V[1]\hookrightarrow C(V)$. 
Just to explain the statement of Theorem~\ref{thmA}
we observe that the condition 
$F\circ \delta=0$ implies $F_{1}\circ d=0$ and then $F_{1}$ induce a 
map in cohomology $\theta\colon H^{*}(V)\to H^{*}(V')=V'$.\\

\bigskip
\section{Proof of Theorem~\ref{thmA}}
Let $X$ be a complex manifold of dimension $n$;  
consider the graded vector space $L=\oplus L^{p}$, where 
$L^{p}=\Gamma(X,\sA^{0,p+1}(T_{X}))$, $-1\le p\le n-1$, and 
two linear maps of degree $+1$, $d\colon L\to L$, $Q\colon 
\odot^{2}L\to L$ defined in the following way: if 
$z_{1},\ldots,z_{n}$ are local holomorphic coordinates, then 
\[ d\left(\phi\desude{~}{z_{i}}\right)=
(\debar\phi)\desude{~}{z_{i}},\qquad \phi\in \sA^{0,*}.\]
If $I,J$ are ordered subsets of $\{1,\ldots,n\}$, 
$a=fd\bar{z}_{I}\desude{~}{z_{i}}$, 
$b=gd\bar{z}_{J}\desude{~}{z_{j}}$, $f,g\in\sA^{0,0}$ then 
\[ Q(a\odot b)=(-1)^{\bar{a}}d\bar{z}_{I}\wedge d\bar{z}_{J}
\left(f\desude{g}{z_{i}}\desude{~}{z_{j}}-
g\desude{f}{z_{j}}\desude{~}{z_{i}}\right),\qquad \bar{a}=\deg(a,L).\]
The equation (\ref{codiff}), with $L$ in place of $V[1]$, gives a 
codifferential $\delta$ of degree 1 on $\bar{S}(L)$ and the 
differential graded coalgebra $(\bar{S}(L),\delta)$ is exactly the 
$L_{\oo}$-algebra associated to the Kodaira-Spencer DGLA $KS_{X}$.\\

If $\Der^{p}(\sA^{*,*},\sA^{*,*})$ denotes the  vector space of 
$\C$-derivations of degree $p$ of the sheaf of graded algebras
$(\sA^{*,*},\wedge)$, where the degree of a $(p,q)$-form is $p+q$
(note that $\de,\debar\in
\Der^{1}(\sA^{*,*},\sA^{*,*})$), then  we can define a morphism of
graded vector spaces
\[ L\mapor{\hat{~}}\Der^{*}(\sA^{*,*},\sA^{*,*})=\somdir{p}{}
\Der^{p}(\sA^{*,*},\sA^{*,*})
,\qquad a\to
\h{a}\] given in local coordinates by 
\[ \h{\phi\desude{~}{z_{i}}}(\eta)=\phi\wedge\left(\desude{~}{z_{i}}
\contr\eta\right).\]
If $\,\bar{a}=p$ then $\h{a}$ is a bihomogeneous 
derivation of bidegree $(-1,p+1)$: in 
particular $\h{a}(\sA^{0,*})=0$.\\

\begin{lem} If $[~,~]$ denotes the standard bracket on 
$\Der^{*}(\sA^{*,*},\sA^{*,*})$, then for every pair of homogeneous 
$a,b\in L$ we have:
\begin{enumerate}
	\item  $\h{da}=[\debar,\h{a}]=\debar\h{a}-(-1)^{\bar{a}}
	\h{a}\debar.$

	\item  $\h{Q(a\odot b)}=-[[\de,\h{a}],\h{b}]=
	(-1)^{\bar{a}}\h{a}\de\h{b}+(-1)^{\bar{a}\,\bar{b}+\bar{b}}\,\h{b}\de\h{a}
	\pm \de\h{a}\h{b}\pm\h{b}\h{a}\de$.
\end{enumerate}
\end{lem}

\begin{pf} By linearity we may assume 
$a=fd\bar{z}_{I}\desude{~}{z_{i}}$, 
$b=gd\bar{z}_{J}\desude{~}{z_{j}}$, $f,g\in\sA^{0,0}$. 
Moreover all the four expressions are derivations vanishing on the 
subalgebra $\sA^{0,*}$ and therefore it is sufficient to check the 
above equalities when computed on the $dz_{i}$'s; 
since $\debar dz_i=\de dz_i=\h{a}\h{b}dz_i=0$,
the computation
becomes  straightforward and it is left to the reader.\qed\end{pf}
~\\
{\em Remark.} The apparent asymmetry in the right hand side of Item 2 
of the above lemma is easily understood: in fact $[\h{a},\h{b}]=0$ and 
then by Jacobi identity 
\[ 0=[\de,[\h{a},\h{b}]]=[[\de,\h{a}],\h{b}]-(-1)^{\bar{a}\,\bar{b}}
[[\de,\h{b}],\h{a}].\]
~\\
Assume now that $X$ is compact K\"ahler, 
fix a K\"ahler metric on $X$ and denote by:
$A^{p,q}=\Gamma(X,\sA^{p,q})$ the vector space of global 
$(p,q)$-forms,
$\debar^{*}\colon A^{p,q}\to A^{p,q-1}$
the adjoint operator of $\debar$, 
$\Delta_{\debar}=\debar\,\debar^{*}+\debar^{*}\debar$ the $\debar$-Laplacian, 
$G_{\debar}$ the associated Green operator, $\sH\subset A^{*,*}$
the graded vector space of harmonic forms, 
$i\colon \sH\to A^{*,*}$ the inclusion and 
$h=Id-\Delta_{\debar} G_{\debar}=Id-G_{\debar}\Delta_{\debar}
\colon A^{*,*}\to \sH$ the harmonic projector.

We identify the graded vector space $M_X$ with the space of 
endomorphisms of harmonic forms $\Hom_{\C}^{*}(\sH,\sH)$. We also 
denote by $N=\Hom_{\C}^{*}(A^{*,*},A^{*,*})$ the graded associative algebra  
of linear endomorphisms of the space of global differential forms on 
$X$.\\ 
For notational simplicity we identify  
$\Der^{*}(\sA^{*,*},\sA^{*,*})$ with its image into $N$.\\

Setting $\tau=G_{\debar}\debar^{*}\de\in N^{0}$ we have  by K\"ahler 
identities (cf. \cite{G-H}, \cite{Weil}):
\[ h\de=\de h=\tau h=h\tau=\de\tau=\tau\de=0\]
\[ [\de,\debar^{*}]=[\de,G_{\debar}]=[\debar,G_{\debar}]=0,\qquad 
[\debar,\tau]=\debar 
G_{\debar}\debar^{*}\de-G_{\debar}\debar^{*}
\de\debar=G_{\debar}\Delta_{\debar}\de=\de.\]

We introduce  the morphism  
\[ F_{1}\colon L\to M_X, \qquad F_{1}(a)=h\h{a}i.\]
We note that $F_{1}$ is a morphism of complexes, in fact 
$F_{1}(da)=h\h{da}i=h(\debar\h{a}\pm\h{a}\debar)i=0$. Next we define, 
for every $m\ge 2$, the morphisms of graded vector spaces 
\[f_{m}\colon \tensor{}{m}L\to M_X,\qquad 
F_{m}\colon \symm{}{m}L\to M_X,\qquad 
F=\sum_{m=1}^{\oo}F_{m}\colon \bar{S}(L)\to M_X,\]
\[ f_{m}(a_{1}\otimes a_{2}\otimes\ldots\otimes a_{m})=
h\h{a_{1}}\tau\h{a_{2}}\tau\h{a_{3}}\ldots\tau\h{a_{m}}i.\]
\[ F_{m}(a_{1}\odot a_{2}\odot\ldots\odot a_{m})=
\sum_{\sigma\in\Sigma_{m}}\epsilon(L,\sigma;a_{1},\ldots,a_{m})
f_{m}(a_{\sigma_{1}}\otimes\ldots\otimes a_{\sigma_{m}}).\]

\begin{thm}\label{1.2} In the above notation $F\circ\delta=0$ and therefore 
\[\Theta=\sum_{m=1}^{\oo}\frac{1}{m!}F^{\odot m}\circ\Delta^{m-1}_{C(KS_X)}
\colon (C(KS_X),\delta)\to (C(M[-1]_X),0)\] 
is an $L_{\oo}$-morphism with linear term $F_{1}$.\end{thm}

\begin{pf} We need to prove that 
for every $m\ge 2$ and $a_{1},\ldots,a_{m}\in L$ we 
have   
\[F_{m}\left(\sum_{\sigma\in S(1,m-1)}\epsilon(L,\sigma)
da_{\sigma_{1}}\odot a_{\sigma_{2}}\odot \ldots\odot 
a_{\sigma_{m}}\right)=\]
\[=-F_{m-1}\left(\sum_{\sigma\in S(2,m-2)}\epsilon(L,\sigma)
Q(a_{\sigma_{1}}\odot a_{\sigma_{2}})\odot a_{\sigma_{3}}\odot \ldots\odot 
a_{\sigma_{m}}\right),\]
where $\epsilon(L,\sigma)=\epsilon(L,\sigma;a_{1},\ldots,a_{m})$.\\
It is convenient to introduce the auxiliary operators
$q\colon \bigotimes^{2}L\to N[1]$, $q(a\otimes 
b)=(-1)^{\bar{a}}\h{a}\de\h{b}$ and 
$g_{m}\colon \bigotimes^{m}L\to M[1]_X$,
\[ g_{m}(a_{1}\otimes \ldots\otimes a_{m})=
-\sum_{i=0}^{m-2}(-1)^{\bar{a_{1}}+\bar{a_{2}}+\ldots+\bar{a_{i}}}
h\h{a_{1}}\tau\ldots\h{a_{i}}\tau q(a_{i+1}\otimes 
a_{i+2})\tau\h{a_{i+3}}\ldots\tau\h{a_{m}}i.\]
Since for every choice of operators
$\alpha=h,\tau$ and $\beta=\tau,i$ and every
$a,b\in L$ we have
\[ \alpha\h{Q(a\odot b)}\beta=\alpha((-1)^{\bar{a}}\h{a}\de\h{b}+
(-1)^{\bar{a}\,\bar{b}+\bar{b}}\h{b}\de\h{a})\beta=
\alpha(q(a\otimes b)+(-1)^{\bar{a}\,\bar{b}}q(b\otimes a))\beta,
\]
a straightforward computation about symmetrization and unshuffles 
gives
\[\sum_{\sigma\in\Sigma_{m}}\epsilon(L,\sigma)
g_{m}(a_{\sigma_{1}}\otimes\ldots\otimes a_{\sigma_{m}})=
-F_{m-1}\left(\sum_{\sigma\in S(2,m-2)}\epsilon(L,\sigma)
Q(a_{\sigma_{1}}\odot a_{\sigma_{2}})\odot a_{\sigma_{3}}\odot \ldots\odot 
a_{\sigma_{m}}\right).\]
On the other hand 
\renewcommand{\arraystretch}{2.7}
\[ \begin{array}{l}
\ds f_{m}\left(\sum_{i=0}^{m-1}(-1)^{\bar{a_{1}}+\ldots+\bar{a_{i}}}
a_{1}\otimes\ldots\otimes a_{i}\otimes 
da_{i+1}\otimes\ldots\otimes a_{m}\right)=\\
\ds =\sum_{i=0}^{m-1}(-1)^{\bar{a_{1}}+\ldots+\bar{a_{i}}}
h\h{a_{1}}\ldots\h{a_{i}}\tau(\debar 
\h{a_{i+1}}-(-1)^{\bar{a_{i+1}}}\h{a_{i+1}}\debar)\tau\ldots \tau 
\h{a_{m}}i\\
\ds =\sum_{i=0}^{m-2}(-1)^{\bar{a_{1}}+\ldots+\bar{a_{i}}}
h\h{a_{1}}\ldots\h{a_{i}}\tau
(-(-1)^{\bar{a_{i+1}}}\h{a_{i+1}}\debar\tau\h{a_{i+2}}+
(-1)^{\bar{a_{i+1}}}\h{a_{i+1}}\tau\debar\h{a_{i+2}})
\tau\ldots \tau 
\h{a_{m}}i\\
\ds =g_{m}(a_{1}\otimes\ldots\otimes a_{m}).\end{array}\]
\renewcommand{\arraystretch}{1}
Taking the symmetrization of this equality we get 
\[ \sum_{\sigma\in\Sigma_{m}}\epsilon(L,\sigma)
g_{m}(a_{\sigma_{1}}\otimes\ldots\otimes a_{\sigma_{m}})=
F_{m}\left(\sum_{\sigma\in S(1,m-1)}\epsilon(L,\sigma)
da_{\sigma_{1}}\odot a_{\sigma_{2}}\odot \ldots\odot 
a_{\sigma_{m}}\right).\]
\qed\end{pf}

Since it is clear that $F_{1}$ is a morphism of complexes inducing the 
morphism $\theta$ in cohomology, Theorem~\ref{thmA} is proved. 
\medskip
~\\
{\em Remark.} If $X$ is a Calabi-Yau manifold with holomorphic volume 
form $\Omega$, then the composition of $F$ with the evaluation at 
$\Omega$ induces an $L_{\oo}$-morphism $C(KS_X)\to C(\sH[n-1])$.\\
For every $m\ge 2$, $\operatorname{ev}_{\Omega}\circ F_{m}\colon 
\bigodot^{m}L\to \sH[n]$ vanishes on 
$\bigodot^{m}\{a\in L\,|\, 
\de(a\contr\Omega)=0\}$.
\medskip
~\\

The following corollary gives a formality criterion:

\begin{cor} In the notation of introduction, 
if $\theta\colon H^{*}(X,T_{X})\to M[-1]_X$ is injective, then $KS_X$ is 
$L_{\oo}$-quasiisomorphic to an abelian differential graded Lie algebra.\end{cor}
	
\begin{pf} Let $H\subset M[-1]_X$ be the image of $\theta$ and let 
$p\colon M[-1]_X\to H$ be a linear projection. 
Since $p$ is a morphism of 
DGLA, the composition $C(KS_X)\mapor{\Theta}C(M[-1]_X)\mapor{p}C(H)$ is 
an $L_{\oo}$ quasiisomorphism.\end{pf}

\bigskip

\section{Applications to deformation theory}

All the technical tools  used in this section are standard and well exposed in the 
literature.
Let $\Art$ be the category of local Artinian $\C$-algebras $(A, 
m_{A})$ with residue field $A/m_{A}=\C$. Following \cite{Sch},
by a functor of Artin rings 
we intend a covariant functor $\sF\colon \Art\to\Set$ such that 
$\sF(\C)=\{0\}$ is a set of cardinality 1.\\
With the term Schlessinger's condition we mean one of the four  conditions 
$(H_{1}),\ldots,(H_{4})$ described in Theorem 2.1 of \cite{Sch}.

\begin{lem}\label{LemmaSch}
Let $\alpha\colon \sF\to \sG$ be a natural transformation of functors 
of Artin rings; if $\sF$ satisfies Schlessinger's conditions $(H_{1})$ and
$(H_{2})$, $\sG$ is prorepresentable and $\alpha\colon t_{\sF}\to t_{\sG}$ 
is injective, then also $\sF$ is prorepresentable.\end{lem}

\begin{pf} 
Since $\sG$ is prorepresentable its tangent space $t_{\sG}$ is finite 
dimensional and then the same holds for  $t_{\sF}$. Moreover for every small 
extension $0\to J\to A\to B\to 0$ there exists a natural 
transitive free action (cf. \cite[2.15]{Sch}) of 
$t_{\sG}\otimes J$ on the nonempty fibres of $\sG(A)\to \sG(B)$. 
Therefore also $t_{\sF}\otimes J$ acts without fixed points on 
$\sF(A)$ and then, according 
to Theorem 2.11 of \cite{Sch}, $\sF$ is prorepresentable.
\qed\end{pf}

For every differential 
graded complex Lie algebra $K=\oplus K^{i}$, we denote respectively by 
$\MC_{K},\Def_{K}\colon \Art\to\Set$ the associated 
Maurer-Cartan and deformation functors  (cf. 
\cite{GoldMil1}, \cite{GoldMil2},  \cite{EDF}):
\[ \MC_{K}(A)=\left\{ a\in K^{1}\otimes m_{A}\, \left|\, 
da+\frac{1}{2}[a,a]=0\,\right.\right\},\qquad 
\Def_{K}(A)=\frac{\MC_{K}(A)}{\exp(K^{0}\otimes m_{A})}.\]
The functors $\MC_{K}$ and $\Def_{K}$ 
are functors of Artin rings satisfying the Schlessinger's conditions 
$(H_{1})$, $(H_{2})$ 
(cf. \cite{Sch},\cite{FM1}), the projection $\MC_{K}\to\Def_{K}$ is 
smooth  and the tangent space $t_{\Def_{K}}$ 
of $\Def_{K}$ is naturally isomorphic to $H^{1}(K)$.\\ 

\begin{example}
\begin{enumerate}
	\item If $K$ has trivial bracket and trivial differential then the gauge 
action is trivial and therefore, for every $(A,m_{A})\in \Art$,  
$\Def_{K}(A)=\MC_{K}(A)=K^{1}\otimes m_{A}$; in particular if 
$K^{1}$ is finite dimensional then $\Def_{K}$ is prorepresented by a smooth germ.
\item  If $K=KS_{X}$ is the Kodaira-Spencer DGLA of a compact complex 
manifold $X$ then $\Def_{K}$ is isomorphic to the functor $\Def_{X}$ of 
infinitesimal deformations of $X$ (cf. \cite{GoldMil2}).
\end{enumerate}
\end{example}	

The functor $\Def_{K}$ has a natural  obstruction theory with 
obstruction space $H^{2}(K)$: this means that for every 
small extension 
\[ e:\quad 0\mapor{}J\mapor{}A\mapor{p}B\mapor{}0\]
in the category $\Art$ it is given an ``obstruction map''
$ob_{e}\colon \Def_{K}\to H^{2}(K)\otimes J$ such that an element 
$b\in \Def_{K}(B)$ lifts to $\Def_{K}(A)$ if and only $ob_{e}(b)=0$. 
Moreover all the obstruction maps behave functorially with respect 
to morphisms of small extensions (cf. e.g. \cite{Artin}, \cite{FM1}).\\
By definition the {\em primary obstruction map} is the obstruction map 
$q_{2}=ob_{\epsilon}\colon H^{1}(K)\to H^{2}(K)$ relative to the small 
extension 
\[ \epsilon:\quad  0\mapor{}\C\mapor{t^{2}}\frac{\C[t]}{(t^{3})}
\mapor{}\frac{\C[t]}{(t^{2})}\mapor{}0.\]
Concretely, if 
$b\in\MC_{K}(B)$ and $a\in K^{1}\otimes m_{A}$ is  a lifting of $b$, then  
by the Jacobi identity $h=da+[a,a]/2\in K^{2}\otimes J$ is a cocycle and its 
cohomology class $ob_{e}(b)=[h]\in H^{2}(K)\otimes J$ does not depend from the 
choice of $a$. It is easy to prove that 
$ob_{e}(b)=0$ if and only if $b$ can be lifted to 
$\MC_{K}(A)$.\\
The map $ob_{e}$ is invariant under the gauge action (this follows from a 
general result  
\cite[7.5]{FM1} but it is also easy to  prove directly) 
and then factors to a map $ob_{e}\colon \Def_{K}(B)\to 
H^{2}(K)\otimes J$. Since the projection  $\MC_{K}\to \Def_{K}$ is  
smooth, we have that the class of $b$ lifts to $\Def_{K}(A)$ if and 
only if $ob_{e}(b)=0$.\\
The obstruction space $O_{K}\subset H^{2}(K)$ is by definition the 
vector space generated by the images  the maps $(Id\otimes 
f)\circ ob_{e}$, where $f\in \Hom_{\C}(J,\C)$ and $e$ ranges over all 
small extension in $\Art$.
\medskip
~\\
{\em Remark.} If the DGLA $K$ is not formal, 
it may happen that the primary obstruction map vanishes but $O_K\not=0$. 
If $O_{K}^{c}\subset O_{K}$ denotes the subspace generated by 
the obstructions coming from all the curvilinear small extensions 
\[ 0\mapor{}\C\mapor{t^{n}}\frac{\C[t]}{(t^{n+1})}
\mapor{}\frac{\C[t]}{(t^{n})}\mapor{}0\]
then, by the (abstract) 
$T^{1}$-lifting theorem \cite{FaMa2}, $\Def_{K}$ is smooth if and only if 
$O_{K}^{c}=0$ but in general 
$O_{K}^{c}\not=O_{K}$ (cf. \cite[5.7]{FM1}).\\
\medskip

Given two differential graded Lie algebras $K,M$,  
every $L_{\oo}$-morphism $\mu\colon C(K)\to 
C(M)$ induces a natural 
transformation $\t{\mu}\colon \Def_{K}\to \Def_{M}$ 
(see e.g. \cite{K}, \cite{EDF}).
Writing $\mu=\sum_{i\le j}\mu^{i}_{j}$, 
$\mu^{i}_{j}\colon \bigodot^{j}K[1]\to\bigodot^{j}M[1]$, 
the morphism $\mu_{1}^{1}$ is a morphism of complexes, 
$H^{1}(\mu_{1}^{1})\colon H^{1}(K)\to H^{1}(M)$ equals the 
restriction of $\t{\mu}$ on tangent spaces and 
$H^{2}(\mu_{1}^{1})\colon H^{2}(K)\to H^{2}(M)$ 
commutes with $\t{\mu}$ and 
all the obstruction maps.  

\begin{prop}\label{vanish} 
Let $K$ be a differential graded Lie algebra, $M\oplus M^{i}$ be a 
graded vector space considered as a differential graded Lie algebra 
with trivial bracket and differential and let
$\mu=\sum_{i\le j}\mu^{i}_{j}\colon C(K)\to C(M)$ be an 
$L_{\oo}$-morphism. Then:
\begin{enumerate}
	\item  If $M^{1}$ is finite dimensional and $H^{1}(\mu_{1}^{1})$ 
	is injective then $\Def_{K}$ is prorepresentable.

	\item  The obstruction space $O_{K}$ is contained in the kernel of 
	$H^{2}(\mu_{1}^{1})\colon H^{2}(K)\to M^{2}$.
\end{enumerate}
\end{prop}
\begin{pf} The first part follows immediately from Lemma~\ref{LemmaSch}.
The second part follows from the fact that all the obstruction maps 
of the functor $\Def_{M}$ are trivial.\qed\end{pf}

If $X$ is a compact K\"ahler manifold we have,
in the notation of the Introduction and Section 2, for every 
$A\in \Art$,
\[ \Def_{X}(A)=\Def_{KS_{X}}(A)=\frac{\ds
\left\{a\in L^{0}\otimes m_{A}\,\left|\, da+\frac{1}{2}Q(a\odot 
a)=0\right.\right\}}{\exp(L^{-1}\otimes m_{A})},\]
\[ \Def_{M[-1]_{X}}(A)=M^{0}_{X}\otimes m_{A}\]
and the natural transformation $\t{F}\colon 
\Def_{KS_{X}}\to\Def_{M[-1]_{X}}$ 
associated  to the $L_{\oo}$-morphism $\Theta$ of Theorem~\ref{1.2} is induced 
by 
\[ \t{\Theta}(a)=\sum_{m=1}^{\oo}\frac{1}{m!}F_{m}(a^{\odot m})=
F(\exp(a)-1), \qquad 
a\in L^{0}\otimes m_{A}.\]

Since $M[-1]_{X}$ carries  the trivial structure of DGLA, 
Proposition~\ref{vanish} gives the 
following result known as the principle ``Obstructions to deformations 
of compact K\"{a}hler manifolds 
annihilate ambient cohomology'' (cf. \cite[10.1]{clemens}, 
\cite[3.5]{ran}).

\begin{cor}\label{principle} Let $X$ be a compact K\"ahler manifold 
and denote by  
$O$ the kernel of 
\[\theta_{2}\colon H^{2}(X,T_{X})\to 
\somdir{r,s}{}\Hom_{\C}^{*}
(H^{r}(\Omega^{s}_{X}),H^{r+2}(\Omega^{s-1}_{X})).\] 
Then for every small extension
$e: 0\mapor{}J\mapor{}A\mapor{p}B\mapor{}0$ 
and  every $b\in\Def_{X}(B)$, the obstruction $ob_{e}(b)$ belongs 
to $O\otimes J$.\end{cor}

\begin{pfof}{Corollary~\ref{CorB}}
We first recall that, if $\sY\to \sB$ is the Kuranishi family of a 
compact complex  manifold $Y$ and $O\subset H^2(Y,T_Y)$ is the subspace generated 
by  all the obstruction to the deformations of $Y$, 
then the singularity $\sB$ is analytically isomorphic to $q^{-1}(0)$, where 
$q\colon (H^1	(Y,T_Y),0)\to (O,0)$  is the Kuranishi map.\\ 
The pull-back of forms and vector fields give a morphism 
of differential graded Lie algebras $\pi^{*}\colon KS_{Y}\to KS_{X}$.
The composition of $\pi^{*}$ with $\Theta$ gives an $L_{\oo}$-morphism 
from $KS_{Y}$ to $M[-1]_X$. It is now sufficient to apply 
Proposition~\ref{vanish}. 
\qed\end{pfof}

\begin{example}
Let $Z$ be a projective Calabi-Yau manifold of dimension $n\ge 3$ 
with $H^{2}(\Oh_{Z})=0$ and let $\pi\colon Y\to Z$ be a smooth Galois double 
cover. Denoting by $D\subset Z$ the branching divisor, $R\subset Y$ 
the ramification divisor and $\pi_{*}\Oh_{Y}=\Oh_{Z}\oplus 
\Oh_{Z}(-L)$ the eigensheaves decomposition we have (cf. \cite{Catacime}, \cite{Par})
$\Oh_{Y}(R)=K_{Y}=\pi^{*}\Oh_{Z}(L)$, $\Oh_Z(D)=\Oh_Z(2L)$, an exact sequence of 
sheaves over $Y$
\[ 0\mapor{}T_{Y}\mapor{}\pi^{*}T_{Z}\mapor{}\Oh_{R}(2R)\mapor{}0\]
and, for every $i$, $H^{i}(\pi^{*}T_{Z})=H^{i}(T_{Z})\oplus 
H^{i}(T_{Z}(-L))$, $H^{i}(\Oh_{R}(2R))=H^{i}(\Oh_{D}(D))$.\\
If $L$ is sufficiently ample then 
$H^{1}(\Oh_{D}(D))=H^{2}(\Oh_{Z})=0$, $H^{2}(T_{Z}(-L))=0$ and then 
$H^{2}(T_{Y})$ injects into $H^{2}(T_{Z})$. Therefore the cup product 
with the pull-back of the holomorphic volume form of $Z$  is nondegenerate and then 
$\theta_{2}\colon H^{2}(T_{Y})\to M^{1}$ is injective. Applying  
Corollary~\ref{CorB} (with $X=Y$) we get that $Y$ is unobstructed.
\end{example}

\bigskip
%\section{biblio}

\begin{tabular}{l} 
Marco Manetti\\
 Dipartimento di Matematica ``G. Castelnuovo'',\\
 Universit\`a di Roma ``La Sapienza'',\\ 
 Piazzale Aldo Moro 5, I-00185 Roma, Italy.\\
 manetti@mat.uniroma1.it,~~~~~
http://www.mat.uniroma1.it/people/manetti/\\
 \end{tabular}

\end{document}